\documentclass[12pt]{article}
\usepackage{graphicx}
\font\tendb=msbm10 at 12pt
\font\sevendb=msbm10 at 9pt
\font\fivedb=msbm10 at 7pt
\newfam\dbfam
\textfont\dbfam=\tendb
\scriptfont\dbfam=\sevendb
\scriptscriptfont\dbfam=\fivedb
\def\db{\fam\dbfam\tendb}
\font\eufm=eufm10\font\eufms=eufm10\font\eufmss=eufm10\newfam\eufam
\textfont\eufam=\eufm\scriptfont\eufam=\eufms\scriptscriptfont\eufam=\eufmss

\font\tendbb=msbm10 at 12pt
\font\sevendbb=msbm7 at 9pt
\font\fivedbb=msbm5 at 6pt
\newfam\dbbfam
\textfont\dbbfam=\tendbb
\scriptfont\dbbfam=\sevendbb
\scriptscriptfont\dbbfam=\fivedbb

 \def \Z {{\db Z}}
 \def \R {\hbox{\db R}}
 
 \def \C {\hbox{\db C}}
 
 \def \S {S^{3}}

 \newcommand{\ann}{I\!\!\!\!F_{p}}
\font\tenMmm=eusm10 at 12pt
\newfam\Mmmfam
\textfont\Mmmfam=\tenMmm

\def\illu #1 by #2 (#3){
  \vbox to #2{
    \hrule width #1 height 0pt depth 0pt
    \vfill
    \special{illustration #3} 
    }
  }

\begin{document}

\null
\vspace{4cm}

\begin{center}  {\large {\bf  A new criterion for  knots with free periods}}\\
 Nafaa Chbili

\end{center}
\vspace{10mm}

 \begin{footnotesize} {\bf R\'esum\'e.} Soient $p\geq 2$ et $q\neq 0$
 deux entiers. Un n\oe ud $K$ de la sph\`ere $\S$ est dit
 $(p,q)$-lenticulaire s'il est invariant par l'action lenticulaire
 $\varphi_{p,q}$. Dans ce travail, nous \'etudions le comportement
 du polyn\^ome de HOMFLY des n\oe uds lenticulaires. Nous
 d\'emontrons que la sym\'etrie lenticulaire  est refl\'et\'ee d'une
 fa\c{c}on tr\`es nette par
 le second coefficient du polyn\^ome de HOMFLY. Comme
 application,  nous d\'emontrons  que 80 parmi les 84 n\oe uds ayant un nombre
 de
 croisements inf\'erieur ou \'egal \`a 9, ne sont pas $(5,1)$-lenticulaires.\\

               {\bf Abstract.} Let $p\geq 2$ and $q\neq
0$ an integer. A knot $K$ in the three-sphere is said to
be a $(p,q)$-lens knot if and only if it covers a link in
the lens space  $L(p,q)$. In this paper, we use the
second coefficient of the HOMFLY polynomial to provide a
necessary  condition for a knot to be a $(p,q)$-lens
knot. As an application, it is shown that this criterion
rules out the possibility of being $(5,1)$-lens for 80
among the 84 knots with less than 9 crossings.

               \end{footnotesize}

\begin{center}{\sc 1- Introduction}\end{center}
This paper is concerned with the question of whether  the symmetry of knots and links in
 the three-sphere is reflected on the quantum invariants. The symmetry we consider in the
  present paper is the free periodicity. A link $L$ in  $\S$ is said to be $p$-freely periodic
   ($p \geq 2$ an integer ) if and only if $L$ is fixed by  an orientation preserving action of
    the finite cyclic group $G=\Z/p\Z$ on the
    three-sphere without  fixed points. It has been
     conjectured since many decades that such an action is topologically conjugate to an
     orthogonal action. Consequently, we are going to limit our interest to links which
     arise as covers of links in the lens space $L(p,q)$. Such a link will be called here
      a $(p,q)$-lens link.\\ The two variable HOMFLY (called also, skein and HOMFLYPT)
polynomial is an invariant of ambiant isotopy of oriented
links, which generalizes both the Alexander and the Jones
polynomials, and can be defined by the following :
$$\begin{array}{ll}
 {\bf (i)}&P_ {\bigcirc}(v,z)=1\\
{\bf (ii)}&v^{-1}P_{L_{+}}(v,z)-vP_{L_{-}}(v,z) =zP_{L_{0}}(v,z),
\end{array}$$
where  $\bigcirc $ is the trivial knot, $L_{+}$,  $L_{-}$ and  $L_{0}$ are three oriented links which are identical except
 near one crossing where they look
like in the following figure:
\begin{center}
\includegraphics[width=8cm,height=2cm]{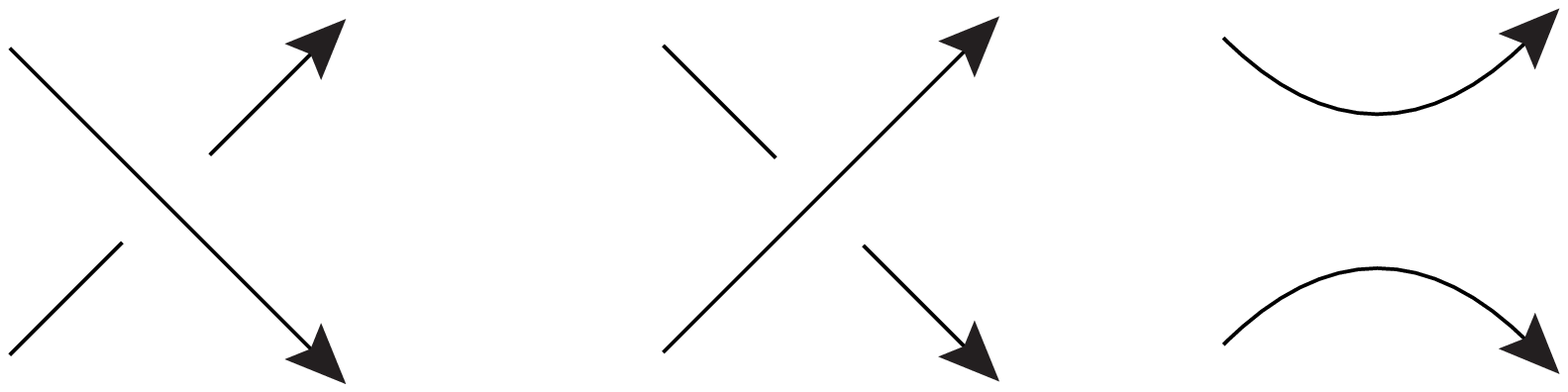}
\end{center}
\begin{center} {\sc Figure 1}
\end{center}

It is well known that the HOMFLY polynomial \cite{LM} takes its values in the ring
$\Z[v^{\pm},z^{\pm}]$. However, if $L$ is a knot then we have $P_{L}(v,z)=
\sum_{i \geq 0}P_{2i,L}(v)z^{2i}$ where $P_{2i,L}$ are elements of $\Z[v^{\pm 2}]$. \\
Knots with free periods were first studied by Hartly \cite{Ha}
who, motivated by a question of R.  Fox, used the Alexander
polynomial to provide a criterion for a knot to be freely
periodic. The first criteria for periodicity of links using the
HOMFLY polynomial is due to Przytycki \cite{Pr}. In \cite{Ch4}, we
used the first term of the HOMFLY polynomial
 to find a necessary condition for a knot to be $p$-freely periodic, for $p$  prime.
  This criterion was applied successfully  to rule out
the possibility of being freely periodic for certain
knots. The aim of this paper is to extend this criterion
to the second coefficient of the HOMFLY polynomial. Thus
we shall prove that similar conditions hold for the
polynomial $P_{2,K}(v)$. The proof of our main result is
based  on the three crucial facts:
\begin{itemize}
\item The combinatorial description of lens knots we provided
in \cite{Ch1}.
\item The techniques developed by Traczyk \cite{Tr1} and Yokota \cite{Yo} in the case of
periodic knots and adapted to freely periodic knots in \cite{Ch4}.
\item The formula for the second term of the HOMFLY
polynomial introduced recently  by Kanenobu and Miyazawa
\cite{KM}.
\end{itemize}
An outline of the present paper is as follows. In section 2 we
introduce our main results. In section 2, basic properties of
freely periodic knots will be summarized. Some properties of the
HOMFLY polynomial, needed in the rest of the paper, are given in
section 4. In section 5, we shall prove Theorem 2.1. In the last
paragraph, our criterion is applied, in the case $p=5$,
to  knots with less than 9 crossings.
\begin{center}{\sc 2- Results and Applications}\end{center}
Let $p$ be a prime and $\ann$ be the cyclic  finite field of $p$ elements. Throughout the rest of our paper we denote by  $P_{2,K}(v)_{p}$ the second
term of the HOMFLY polynomial considered with coefficients reduced modulo $p$. If $m$ and $n$
are two integers then $T(n,m)$  denotes the torus link of type $(n,m)$. Recall here that the
 number of components of  $T(n,m)$ is equal to gcd$(n,m)$. In particular
 $T(n,m)$ is a knot if and only if $n$ and $m$ are coprime.\\

 {\bf Theorem 2.1.} {\sl Let $p>3$ be a prime, $q=\pm 1$
 and $K$ a $(p,q)$-lens knot. Then $P_{2,K}(v)_{p} \in
\Gamma_{p,q}$, where $\Gamma_{p,q}$ is the $\ann [v^{\pm
2p}]$-module generated by $P_{2,T(\alpha,\alpha q \pm p
)}(v)_{p}$ for all $1\leq \alpha \leq p-1$.}\\

 This result is more significant for small values of $p$.
Indeed, for such values the generators of the module
$\Gamma_{p,q}$ are easily computed using the formula
given by V. Jones \cite{Jo} for the HOMFLY polynomial of
torus knots. This fact is illustrated by the following
corollary:\\
 {\bf Corollary 2.2.} {\sl Let $q=\pm 1$ and  $K$ be a
$(5,q)$-lens knot. Then $P_{2,K}(v)_{5} \in
{I\!\!\!\!F_{5}}[v^{\pm 10}]$-module generated by
$v^{q8}$.}\\

Proof of Corollary 2.2. According to theorem 2.1, the
generators of $\Gamma_{5,1}$ are given by
$P_{2,T(\alpha,\alpha  \pm 5 )}(v)_{5}$ for $1\leq \alpha
\leq 4$. We use the formula given in section 4 to compute
the HOMFLY polynomial of torus knots. These generators
are given by the list below :
$$
\begin{array}{ll}
P_{2,T(1,6)}(v)_{5}=1,&P_{2,T(1,-4)}(v)_{5}=1,\\
P_{2,T(2,7)}(v)_{5}=10v^{6}-4v^{8},&P_{2,T(2,-3)}(v)_{5}=v^{-2},\\
P_{2,T(3,8)}(v)_{5}=105v^{14}-21v^{8}-105v^{6},&P_{2,T(3,-2)}(v)_{5}=v^{-2},\\
P_{2,T(4,9)}(v)_{5}=770v^{24}-1210v^{26}-70v^{30}+56v^{28},&P_{2,T(4,-1)}(v)_{5}=1.
\end{array}
$$
A similar computation can be easily made in the case
$q=-1$.\\
 {\bf Remark 2.3.}  In the case $p=7$, the module
$\Gamma_{7,1}$ is generated by the two elements:
$2v^{6}+3v^{8}$ and $6v^{8}+4v^{10}$. Consequently the
module $\Gamma_{7,-1}$ is generated by $2v^{-6}+3v^{-8}$
and $6v^{-8}+4v^{-10}$.\\
{\bf Application.} Corollary
2.2 provides  a criterion for a knot of the three-sphere
to be fixed by the lens transformation $\varphi_{5,\pm
1}$. Hence, given a knot $K$, if the polynomial
$P_{2,K}(v)_{5}$ does not satisfy the condition given by
corollary 2.2 then $K$ is not a $(5,\pm 1)$-lens knot.
Let us illustrate this by considering the knot
$K=8_{13}$.  According to the table in \cite{LM} we have
$$P_{2,8_{13}}(v)_{5} = v^{-2}-1-2v^{2}+v^{4}.$$ As
$P_{2,8_{13}}(v)_{5}$ is not in the ${I\!\!\!\!F_{5}}[v^{\pm
10}]$-module generated by $v^{8}$. Then $K=8_{10}$ is not a
$(5,1)$-lens knot. It is worth mentioning that the criterion we
introduced in \cite{Ch4} using the first coefficient of the HOMFLY
polynomial does not decide in the case of the knot $K=8_{13}$.
Thus, corollary 2.2 is not a consequence of the results we
introduced in \cite{Ch4}. More applications are given in the last
section of this
paper.\\
\begin{center}{\sc 3- Freely periodic links}\end{center}
Symmetry of knots and links is a vast subject that has
fascinated researchers since the early age of knot
theory. Problems as chirality and invertibility have
motivated classical knot theory for a long time. Roughly,
a knot in $\S$ is said to be symmetric if and only if $K$
is fixed by an action of a finite cyclic group on $\S$.
According to the set of fixed points of the action, we
can distinguish many kinds of symmetry. In this section
we focus on the case where the action has no fixed
points. We define freely periodic knots then we review
some basic properties of this family of knots and
links.\\
 {\bf Definition 3.1.} Let $p\geq 2$ be an
integer. A link $L$ in $\S$ is said to be $p$-freely
periodic if and only if there exists an orientation
preserving diffeomorphism $h:\S \longrightarrow \S$ such
that:\\
1) $h^{i}$ has no fixed points for all $1\leq i
\leq p-1$,\\
2) $h^{p}=Id_{\S},\\
$ 3) $h(L)=L$.\\
 {\bf Example 3.2.}
Let $L$ be the torus knot $T(2,5)$. It is well known that
$L$ can  be seen as the intersection between an
appropriate three-sphere and the complex surface defined
by :
$$\Sigma=\{(z_{1},z_{2}) \in \C
\times \C ; z_{1}^{2}+z_{2}^{5}=0\}.$$ Let us consider
the diffeomorphism :
$$\begin{array}{cccl}
h:& S^{3} & \longrightarrow & S^{3} \\
   & (z_{1},z_{2}) & \longmapsto & (e^{\frac{2i\pi}{3}}z_{1},
e^{\frac{2i\pi}{3}}z_{2}).
\end{array} $$
Obviously $h$ satisfies conditions 1 and 2 of  definition 3.1. Moreover, one may easily check that $h(L)=L$. Thus $L$ is a freely periodic knot with period 3.\\
{\bf Remark 3.3.} Let $p \geq 2$ and $q$ an  integer
such that gcd$(p,q)=1$. Consider $\varphi_{p,q}$ the diffeomorphism given by:
$$\begin{array}{cccl}
\varphi_{p,q} :& S^{3} & \longrightarrow & S^{3} \\
   & (z_{1},z_{2}) & \longmapsto & (e^{\frac{2i\pi}{p}}z_{1},
e^{\frac{2iq\pi}{p}}z_{2}).
\end{array} $$
It is easy to see that $\varphi_{p,q}$ is an orientation
preserving diffeomorphism of order $p$ and  that
$\varphi_{p,q}$ has no fixed point. Moreover, we have a
$p$-fold cyclic covering $(\pi_{p,q},\S, L(p,q)$).\\ {\bf
Definition 3.4.} {\sl  Let $p \geq 2$ and $q$ an integer
such that gcd$(p,q)=1$. A link $L$ of $\S$ is said to be
a $(p,q)$-lens link if and only if $L$ is mapped onto
itself by $\varphi_{p,q}$.}\\ It is worth mentioning that
lens links are the only examples we know of freely
periodic links. More precisely we have the following
conjecture proved for $p=2$ and 3.\\ {\bf Conjecture 3.5
\cite{Ru}.} {\sl Let $p$ be a prime and $h:S^{3}
\longrightarrow  S^{3}$ an orientation preserving
diffeomorphim of order $p$ such that  for all $1\leq i
\leq p-1$, $h^{i}$ has no fixed points. Then there exists
an integer $q$ such that $h$ is topologically conjugate
to $\varphi_{p,q}$.} \\ Let $n\geq 1$ be an integer. An
$n$-tangle $T$ is a submanifold of dimension one
 in $\R^{2}\times I$ such that the boundary of $T$ is made up of $2n$ points
  $\{A_{1},\dots,A_{n}  \}\times\{0,1\}$. If  $T$ and $T'$ are two $n$-tangles we define the
   product $TT'$ by putting $T$ over $T'$ as follows:
\vspace{1cm}
$$\begin{picture}(0,0)
\put(-20,0){\framebox(40,20){\footnotesize {T}}}
\put(-20,-20){\framebox(40,20){\footnotesize {T'}}}
\end{picture}$$\\

As in the case of braids we define the closure of $T$ and
we denote by $\widehat T$ the link obtained from $T$ by
joining $A_{i}\times {1}$ to $A_{i}\times {0}$ by a
simple arc without adding any crossing. Throughout the
rest of this paper ${\cal B}_{n}$ denotes the $n$-string
braid group. It is well known that this group has the
following presentation:
$$B_{n}= \langle \sigma_{1}, \sigma_{2}, \dots ,\sigma_{n-1} | \sigma_{i} \sigma_{j}=\sigma_{j} \sigma_{i} \mbox{ if } |i-j| \geq 2 {\mbox{ and }}
\sigma_{i} \sigma_{i+1} \sigma_{i}=\sigma_{i+1} \sigma_{i}\sigma_{i+1},
\forall 1 \leq i \leq n-2
\rangle. $$
For $n>2$, the group ${\cal B}_{n}$ is not abelian. Its
center is known to be generated  by the element
$\Omega_{n} = (\sigma_{1} \sigma_{2} ...
\sigma_{n-1})^{n}$. The following theorem gives a
combinatorial description of lens links.\\
{\bf Theorem
3.6. \cite{Ch1}}  {\sl A link $K$ of $\S$ is a
$(p,q)$-lens link if and only if there exists an integer
$n \neq 0$ and an $n$-tangle $T$ such that}:
$$K=\widehat{T^{p}(\sigma_{1} \sigma_{2}...\sigma_{n-1})^{nq}}.$$
\\
\vspace{3cm}
\begin{picture}(0,0)
\put(210,0){\framebox(40,20){\footnotesize {T}}}
\put(210,-20){\framebox(40,20){\footnotesize {T}}}
\put(230,-25){.}
\put(230,-30){.}
\put(230,-35){.}
\put(210,-60){\framebox(40,20){\footnotesize {T}}}
\put(210,-80){\framebox(40,20){\footnotesize {$\Omega_{n}^{q}$}}}
\end{picture}

{\bf Remark 3.7.} Let $n$ and $m$ be two integers. The
torus link $T(n,m)$ is the closure of the braid
$(\sigma_{1} \sigma_{2} ... \sigma_{n-1})^{m}$. Using
elementary techniques we can prove that $T(n,m)$ is a
$(p,q)$-lens link if and only if $p$ divides $m-nq$.\\
\begin{center}{\sc 4- The HOMFLY polynomial}\end{center}
The discovery of the Jones polynomial \cite{Jo} led to a
significant progress in knot theory. The Jones polynomial
was followed by a
 family of invariants of knots and three-manifolds called the quantum invariants.
 Among this family of invariants  the HOMFLY polynomial which is an invariant of ambiant isotopy of oriented
 links. This invariant is a two-variable Laurent polynomial which can be seen as a
  generalization of the Jones and the Alexander
polynomial. This section is to review some properties of
this polynomial needed in the sequel. At the beginning
let us fix some notations.\\
 Let $L=l_{1}\cup l_{2}\cup\dots\cup l_{n}$ be an
$n$-component link of the three-sphere. Throughout the
rest of this paper $\lambda_{i,j}$ denotes the linking
number of the two components $l_{i}$ and $l_{j}$ and
$\lambda$  denotes the total linking number of the link
$L$. It is well known that the HOMFLY polynomials takes
values in the ring $\Z[v^{\pm 1},z^{\pm 1} ]$. Moreover,
we can  write $P_{L}(v,z)=\sum_{i\geq 0}
P_{1-n+2i,L}(v)z^{1-n+2i}$ where $P_{1-n+2i,L} \in
\Z[v^{\pm 2}]$ if $n$ is odd and $P_{1-n+2i,L} \in
\Z[v^{\pm 1}]$ if $n$ is even.\\
{\bf Proposition 4.1 \cite{LM}.} {\sl Let $L=l_{1}\cup
l_{2}\cup\dots\cup l_{n}$ be an $n$-component link then:
$$P_{1-n,L}(v)=v^{2\lambda}(v^{-1}-v)^{n-1}\displaystyle\prod_{i=1}^{n}P_{0,l_{i}}(v).$$}

Motivated by this proposition, Kanenobu and Miyazawa
\cite{KM} introduced a similar formula  for the
polynomial $P_{3-n,L}$.\\ {\bf Theorem 4.2 \cite{KM}.}
{\sl Let $n\geq 3$ be  an integer and  $L=l_{1}\cup
l_{2}\cup\dots\cup l_{n}$ an $n$-component link then:
$$\begin{array}{rl}
P_{3-n,L}(v)=&v^{2\lambda}(v^{-1}-v)^{n-2}\displaystyle\sum_{i<j}(v^{-2\lambda_{i,j}}
P_{1,L_{i,j}}(v)\displaystyle\prod_{k\neq i,j}P_{0,L_{k}}(v))\\
&-(n-2)v^{2\lambda}(v^{-1}-v)^{n-1}\displaystyle\sum_{i=1}^{n}(P_{2,l_{i}}(v)
\displaystyle\prod_{j\neq i}P_{0,l_{j}}(v)),
\end{array}
$$
where $L_{i,j}$ denotes the 2-component link $l_{i}\cup
l_{j}$.}\\ The HOMFLY polynomial of torus knots was
computed by V. Jones. To introduce the Jones  formula, we find
it more convenient to use the polynomial  $X_{L}(q,t)$.
This is a version of the HOMFLY polynomial related to
 $P_{L}(v,z)$ by the variable changes: $z=q^{1/2}-q^{-1/2}$ and  $v=(t q)^{1/2}$.
 Let $k \geq 1$ be an integer, we define:
 $[k]!=(1-q)(1-q^{2})...(1-q^{k})$ and
$[\bar k]=\displaystyle\frac {1-q^{k}}{1-q}.$ \\

{\bf Theorem 4.3 \cite{Jo}.} {\sl For the torus knot
$T(n,m)$ we have:
$$ X_{
T(n,m)}(q,t)=\frac {t ^{(n-1)(m-1)/2}}{[ \bar
n](1-t q)}
\displaystyle \sum_{\stackrel{\normalsize{\gamma + \beta
+1=n}}{\gamma
\geq 0,
\beta \geq 0 }}(-1)^{ \beta }
\displaystyle\frac{q^{\frac {2\beta m +\gamma (\gamma +1)}{2}}}{[\gamma]![\beta]!}
\displaystyle \prod_{i=-\gamma}^{\beta}(q^{i}-\lambda q).$$}
\begin{center}{\sc 5- Proof of Theorem 2.1}\end{center}
Most of the techniques used in this section were first
developed by Traczyk \cite{Tr1} to study the HOMFLY
polynomial of periodic knots (in some sens this class of
knots corresponds to the $(p,0)$-lens knots).  In this
section we aim to adapt these techniques to the case of
freely periodic knots. This will be done in two steps. In
the first one, we prove that $P_{D}$ belongs to $\Gamma_{p,q}'$, where $\Gamma_{p,q}'$
is the $\ann [v^{\pm
2p}]$-module
generated by the polynomials of torus knots $T(n,nq+p)$.
The second step explains how to extract a finite set of
generators for $\Gamma_{p,q}'$. Let us fix some notations.
By $T_{+}$, $T_{-}$ and $T_{0}$, we denote three tangles
which  are identical except near one crossing where they
look like in figure 1. By $D_{+}$ (respectively $D_{-}$,
$D_{0}$)  we denote a diagram of the $(p,q)$-lens link
$D_{+}=\widehat {T_{+}^{p}\Omega_{n}^{q}}$ (respectively
$D_{-}=\widehat { T_{-}^{p}\Omega_{n}^{q}}$ and
$D_{0}=\widehat { T_{0}^{p}\Omega_{n}^{q})}$). It is worth mentioning
that if $D_{+}$ is a knot then $D_{-}$ is also a knot.
However, $D_{0}$ is a link with 2 or $p+1$ components. In
the case $D_{0}$ has $p+1$ components $D_{1}\cup
D_{2}\cup\dots  \cup D_{p+1}$, then one component (say
$D_{1}$)  is invariant by $\varphi_{p,q}$, the others are
cyclically permuted by $\varphi_{p,q}$. We shall prove by
induction on the number of crossings of $D$, that
$P_{2,D} \in \Gamma_{p,q}'$. Let $D$ be a $(p,q)$-lens
diagram. Assume that for all $(p,q)$-lens diagram $D'$
with less crossings than $D$ we have $P_{2,D'} \in
\Gamma_{p,q}$.
In \cite{Ch4}, the following lemma was proved:\\
{\bf
Lemma 5.1.} {\sl  Let $p$ be a prime. The following
congruence holds modulo $p$:
$$ v^{-p}P_{D_{+}}(v,z)-v^{p}P_{D_{-}}(v,z)\equiv z^{p}P_{D_{0}}(v,z). $$}

{\bf Proposition  5.2.} {\sl Let $p \geq 5$ be a prime.\\ i) If
$D_{0}$ has two components then:\\
$v^{-p}P_{2,D_{+}}(v)_{p}-v^{p}P_{2,D_{-}}(v)_{p}=0 $\\
ii) If $D_{0}$ has $p+1$ components then:\\
$v^{-p}P_{2,D_{+}}(v)_{p}-v^{p}P_{2,D_{-}}(v)_{p}=v^{2\lambda}(v^{-1}-v)^{p}
 P_{2,D_{1}}(v)(P_{0,D_{2}}(v))^{p} $}.\\

Proof: According to lemma 5.1 we have the following
congruence modulo $p$:
$$v^{-p}P_{2,D_{+}}(v)_{p}-v^{p}P_{2,D_{-}}(v)_{p}=P_{3-(p+1),D_{0}}(v).$$
Obviously, if $D_{0}$ has two components then $P_{3-(p+1),D_{0}}$
is zero. Assume now  that   $D_{0}$ has $p+1$ components $D_{1}$,$D_{2}$,\dots,$
D_{p+1}$. Let we denote by $D_{i,j}$ the two-component link $D_{i}\cup D_{j}$
and define $G$ and $H$ as follows:
$$\begin{array}{ll}
G(v)=&v^{2\lambda}(v^{-1}-v)^{p-1}\displaystyle\sum_{i<j}(v^{-2\lambda_{i,j}}
P_{1,D_{i,j}}(v)\displaystyle\prod_{k\neq
i,j}P_{0,D_{k}}(v)),\\
H(v)=&-(p-1)v^{2\lambda}(v^{-1}-v)^{p}\displaystyle\sum_{i=1}^{p+1}(P_{2,D_{i}}(v)
\displaystyle\prod_{j\neq i}P_{0,D_{j}}(v)).
\end{array}$$
It can be easily seen from Theorem 4.2 that :
$$v^{-p}P_{2,D_{+}}(v)_{p}-v^{p}P_{2,D_{-}}(v)_{p} \equiv G(v)+ H(v).$$
By the fact that components $D_{2}$,$D_{3}$,\dots,$D_{p+1}$ are
identical and cyclically permuted by the action of $\Z/p\Z$ we can
write:
$$ G(v)=v^{2\lambda}(v^{-1}-v)^{p-1}
(\displaystyle\sum_{j=2}^{p+1}v^{-2\lambda_{1,j}}
P_{1,D_{1,j}}(v)\displaystyle\prod_{k=2,k\neq j}^{p+1}P_{0,D_{k}}(v)+
\displaystyle\sum_{1<i<j}(v^{-2\lambda_{i,j}}
P_{1,D_{i,j}}(v)\displaystyle\prod_{k\neq i,j}P_{0,D_{k}}(v)).
$$
Using the fact that $\lambda_{1,2}=\lambda_{1,j}$ and that $D_{1,2}=D_{1,j}$
for all $2\leq j \leq p+1$, we get:

$$\displaystyle\sum_{j=2}^{p+1}v^{-2\lambda_{1,j}}
P_{1,D_{1,j}}(v)\displaystyle\prod_{k=2,k\neq j}^{p+1}P_{0,D_{k}}(v)=pv^{-2\lambda_{1,2}}
P_{1,D_{1,2}}(v)(P_{0,D_{2}}(v))^{p-1}.
$$
On the other hand :
$$
\begin{array}{ll}
\displaystyle\sum_{1<i<j}(v^{-2\lambda_{i,j}}
P_{1,D_{i,j}}(v)\displaystyle\prod_{k\neq i,j}P_{0,D_{k}}(v))&=
\displaystyle\sum_{1<i<j}(v^{-2\lambda_{i,j}}
P_{1,D_{i,j}}(v)P_{0,D_{1}}(v)(P_{0,D_{2}}(v))^{p-2}\\
&=P_{0,D_{1}}(v)(P_{0,D_{2}}(v))^{p-2}\displaystyle\sum_{1<i<j}v^{-2\lambda_{i,j}}
P_{1,D_{i,j}}(v).
\end{array}$$

One may check easily that $\displaystyle\sum_{1<i<j}v^{-2\lambda_{i,j}}
P_{1,D_{i,j}}(v) \equiv 0$ modulo $p$. Thus G(v) is zero modulo
$p$. A similar computation shows that:
$$
\begin{array}{ll}
\displaystyle\sum_{i=1}^{p+1}(P_{2,D_{i}}(v)
\displaystyle\prod_{j\neq
i}P_{0,D_{j}}(v))&=P_{2,D_{1}}(v)(P_{0,D_{2}}(v))^{p}+
\displaystyle\sum_{i=2}^{p+1}
P_{2,D_{2}}(v)P_{0,D_{1}}(v)(P_{0,D_{2}}(v))^{p-2}\\
&\equiv P_{2,D_{1}}(v)(P_{0,D_{2}}(v))^{p} \mbox {\,\, mod } p.
\end{array}
$$

Therefore: $H(v)\equiv v^{2\lambda}(v^{-1}-v)^{p}
 P_{2,D_{1}}(v)(P_{0,D_{2}}(v))^{p}$ modulo $p$.
 This ends the proof of Proposition 5.2.\\

{\bf Lemma 5.3.} {\sl $P_{2,D_{+}} \in \Gamma_{p,q}'$ if
and only if $P_{2,D_{-}} \in \Gamma_{p,q}'$}.\\

Proof: It is easy to see that the result is true in the case $D_{0}$
has two components. If $D_{0}$ has $p+1$ components $D_{1}$,
$D_{2}$, ..., $D_{p+1}$, then $D_{1}$ is a $(p,q)$-lens diagram with
less crossings  than $D_{0}$. Thus $P_{2,D_{1}} \in \Gamma_{p,q}'$
by the induction assumption. Moreover, an easy computation shows that
$p$ divides the total linking number $\lambda$. According to
Proposition 5.2 we have:
$$
v^{-p}P_{2,D_{+}}(v)_{p}-v^{p}P_{2,D_{-}}(v)_{p}=v^{2\lambda}(v^{-1}-v)^{p}
 P_{2,D_{1}}(v)(P_{0,D_{2}}(v))^{p}.
$$
Obviously, $
v^{2\lambda}(v^{-1}-v)^{p}(P_{0,D_{2}}(v))^{p}$ belongs
to $\Gamma_{p,q}'$. Therefore, the second term in the
previous identity belongs to the module $\Gamma_{p,q}'$.
Consequently, $P_{2,D_{+}} \in \Gamma_{p,q}'$ if and only
if $P_{2,D_{-}} \in \Gamma_{p,q}'$.\\
{\sl Notation.}
Throughout the rest of this paper, we denote by
$D_{+}\leftrightarrow D_{-}$ the operation that consists
of modifying $p$-crossings to transform the diagram
$D_{+}$ into the diagram $D_{-}$ or vice-versa.

The following two lemmas, explain how to use the operation $D_{+}\leftrightarrow D_{-}$
to transform a lens diagram $D$ into  a torus knot. Details
about these techniques can be found in \cite{Ch4}.\\

{\bf Lemma 5.4.} {\sl Every $(p,q)$-lens diagram may be
transformed into a $(p,q)$-lens closed braid by a series of
operations $D_{+}\leftrightarrow D_{-}$ without increasing the
number of crossings.
}\\
{\bf Lemma 5.5.} {\sl Let $B$ be an $n$-braid. The $(p,q)-$lens
braid $B^{p}\Omega_{n}^{q}$ may be transformed into the torus
knot $T(n,nq+p)$ by a series of operations $D_{+}\leftrightarrow
D_{-}$.
}\\
It remains now to explain how to extract a finite set of
generators for the module $\Gamma_{p,q}'$. Our approach here  will
be based on some combinatorial elementary properties of torus
knots. Namely, we shall adapt the $D_{+}\leftrightarrow D_{-}$
operation to diagrams of torus knots. Therefore, an easy induction
will end the proof of Theorem 2.1. We refer the reader to
\cite{Ch4} for more
details.\\
\begin{center}{\sc 6- More applications}\end{center}
This section is devoted to some applications of Theorem 2.1. As
explained earlier in this paper, we can use the criterion Theorem
2.1 provides to decide if a knot $K$ is not a $(p,q)$-lens knot.
Let us first recall that in \cite{Ch4}, we introduced a criterion
for free periodicity using the first coefficient of the HOMFLY
polynomial $P_{0,K}$ (this criterion will be called the
$P_{0}$-criterion). In the case $p=5$, this criterion
writes as follows:\\

{\bf The $P_{0}$-criterion}. {\sl If $K$ is a (5,1)-lens knot,
then $P_{K}=\sum a_{2i}v^{2i}$ with
$a_{10k+4}=2a_{10k+2}$ and $a_{10k+6}=2a_{10k+8}$ for all
$k\in \Z$.}
\\
 Our aim here is to understand how powerful is the
 criterion introduced in section 2 in detecting free
periodicity. Namely, we shall compare the condition
obtained in the present paper to the $P_{0}$-criterion.
To do, let us apply both of them to the 84 knots with less
than 9 crossings. This is explained in the following
table where:\\
 The first column gives the prime knot
according to the notations used in \cite{Ro}.\\ The
second column (resp. the third) provides informations
about the $P_{0}$-criterion (resp. the $P_{2}$-criterion introduced by corollary 2.2) as
follows:\\
D means that the criterion decides that the
knot is not a (5,1)-lens knot.\\
ND means that the
criterion does not decide that the knot is not a (5,1)-lens
knot.
\\
\begin{footnotesize}
$$\begin{array}{llllll}
Knot&P_{0}-\mbox{criterion }&P_{2}-\mbox{criterion
}&Knot&P_{0}-\mbox{criterion } &P_{2}-\mbox{criterion }\\

3_{1}&ND&D&9_{8}&D&D\\ 4_{1}&D&ND&9_{9}&ND&D\\
5_{1}&D&D&9_{10}&ND&D\\5_{2}&D&D&9_{11}&D&D\\
6_{1}&D&D&9_{12}&D&D\\6_{2}&ND&D&9_{13}&ND&D\\
6_{3}&D&D&9_{14}&D&D\\7_{1}&ND&ND&9_{15}&D&D\\
7_{2}&D&D&9_{16}&ND&D\\7_{3}&D&D&9_{17}&D&D\\
7_{4}&D&D&9_{18}&D&D\\7_{5}&D&D&9_{19}&D&D\\
7_{6}&D&D&9_{20}&D&D\\7_{7}&ND&D&9_{21}&D&D\\
8_{1}&D&D&9_{22}&D&D\\8_{2}&D&D&9_{23}&D&D\\

8_{3}&D&D&9_{24}&D&D\\8_{4}&D&D&9_{25}&ND&D\\
8_{5}&D&D&9_{26}&D&D\\8_{6}&D&D&9_{27}&D&D\\
\end{array}$$
$$\begin{array}{llllll}
Knot&P_{0}-\mbox{criterion }&P_{2}-\mbox{criterion
}&Knot&P_{0}-\mbox{criterion } &P_{2}-\mbox{criterion }\\

8_{7}&D&D&9_{28}&D&D\\8_{8}&D&D&9_{29}&D&D\\
8_{9}&D&D&9_{30}&D&D\\8_{10}&D&D&9_{31}&ND&D\\
8_{11}&D&D&9_{32}&D&D\\8_{12}&D&D&9_{33}&ND&D\\
8_{13}&ND&D&9_{34}&D&D\\8_{14}&ND&D&9_{35}&ND&D\\
8_{15}&D&D&9_{36}&D&D\\8_{16}&ND&D&9_{37}&D&D\\
8_{17}&D&D&9_{38}&ND&D\\8_{18}&D&D&9_{39}&D&D\\

8_{19}&ND&ND&9_{40}&ND&D\\8_{20}&ND&D&9_{41}&D&D\\
8_{21}&D&D&9_{42}&D&D\\9_{1}&ND&ND&9_{43}&ND&D\\
 9_{2}&D&D&9_{44}&D&D\\9_{3}&ND&D&9_{45}&D&D\\
9_{4}&D&D&9_{46}&D&D\\9_{5}&D&D&9_{47}&D&D\\
9_{6}&ND&D&9_{48}&D&D\\9_{7}&D&D&9_{49}&ND&D
\end{array}$$
\end{footnotesize}


\begin{footnotesize}
               Nafaa Chbili\\
               D\'epartement de Math\'ematiques,
               Facult\'e des Sciences de Monastir.\\
               Boulevard de l'environnement,
               Monastir 5000, Tunisia.\\
               e-mail {\underline {nafaa.chbili@esstt.rnu.tn}} \end{footnotesize}

\end{document}